\documentclass[11pt]{article}
\input{epsf.sty}
\newtheorem{theorem}{Theorem}[section]

\newtheorem{remark}{Remark}[section]
\newtheorem{example}{Example}[section]

\newtheorem{proposition}{Proposition}[section]

\newcommand{\Prob}{{\rm P}}

\newcommand{\ds}{\displaystyle}

\newcommand{\Stati}{\hbox{$\cal S$}}

\newcommand{\Interi}{\hbox{\bf Z}}

\newcommand{\qed}{\hfill\rule{2mm}{2mm}}      
\newcommand{\remove}[1]{}
\def\eq#1{~{\rm (\ref{equation:#1})}}

\newenvironment{proof}{\begin{trivlist}
\item[\hspace{\labelsep}{\bf\noindent Proof. }]}
{\qed\end{trivlist}}
\title{\Large\bf ON FIRST-PASSAGE-TIME DENSITIES
FOR CERTAIN SYMMETRIC MARKOV CHAINS\footnote{This work has been 
partially supported by MIUR (cofin 2003) and INdAM (G.N.I.M.)}
}
\author{
{\sc A. Di Crescenzo\footnote{Corresponding author; tel.: +39-089-963349; fax: +39-089-963303}}\\
\normalsize Dipartimento di Matematica e Informatica\\ 
\normalsize Universit\`a degli Studi di Salerno\\ 
\normalsize 84084 Fisciano (SA), Italy\\ 
\normalsize Email: {\tt adicrescenzo@unisa.it}
\and
{\sc A. Nastro}\\
\normalsize Dipartimento di Matematica e Applicazioni\\ 
\normalsize Universit\`a di Napoli Federico II\\ 
\normalsize Via Cintia, 80126 Napoli, Italy\\
\normalsize Email: {\tt patrizia.nastro@dma.unina.it}  
}
\begin{document}

\baselineskip 14pt 

\maketitle

\begin{abstract}
The spatial symmetry property of truncated birth-death processes studied 
in Di~Cre\-scenzo \cite{Di98} is extended to a wider family of continuous-time 
Markov chains. We show that it yields simple expressions for first-passage-time 
densities and avoiding transition probabilities, and apply it to a bilateral 
birth-death process with jumps. It is finally proved that this symmetry property 
is preserved within the family of strongly similar Markov chains. 



\end{abstract}
%
\section{Introduction}
A spatial symmetry for the transition probabilities of truncated birth-death processes 
has been studied in Di~Crescen\-zo \cite{Di98}. Such a property leads to 
simple expressions for certain first-passage-time densities and avoiding transition 
probabilities. In this paper we aim to extend those results to a 
wider class of continuous-time Markov chains. 
\par
Given a set $\{x_n\}$ of positive real numbers and the transition probabilities 
$p_{k,n}(t)$ of a continuous-time Markov chain whose state-space is 
$\{0,1,\ldots,N\}$ or $\Interi$, in Section 2 we introduce the following spatial symmetry 
property:   
\begin{equation}
 p_{N-k,N-n}(t)={x_n\over x_k}\,p_{k,n}(t). 
 \label{equation:2}
\end{equation}
In section 3 we point out some properties of first-passage-time densities 
and avoiding transition probabilities for Markov chains that are symmetric in the sense 
of\eq{2}. These properties allow one to obtain simple expressions for first-passage-time 
densities in terms of probability current functions, and for avoiding transition 
probabilities in terms of the `free' transition probabilities. In Section 4 we then apply 
these results to a special bilateral birth-death process with jumps. 
Finally, in Section 5 we refer to the notion of strong similarity 
between the transition probabilities of Markov chains, expressed by 
$\widetilde  p_{k,n}(t)=(\beta_n/\beta_k)\,p_{k,n}(t)$ (see Pollett \cite{Po01}, and 
references therein) and show the following preservation result: if $p_{k,n}(t)$ possesses 
the symmetry property\eq{2}, then also $\widetilde  p_{k,n}(t)$ does it. 
\section{Symmetric Markov chains}
Let $\{X(t),\,t\geq 0\}$ be a homogeneous continuous-time Markov chain on a 
state-space $\Stati$. We shall assume that $\Stati=\{0,1,\ldots,N\}$, where $N$ is 
a fixed positive integer, or $\Stati=\Interi\equiv \{\ldots,-1,0,1,\ldots\}$. Let  
\begin{equation}
	p_{k,n}(t)=\Pr\{X(\tau+t)=n\,|\,X(\tau)=k\}, \qquad 
 k,n\in\Stati; \;\; t,\tau\geq 0
 \label{equation:1}
\end{equation}
be the stationary transition probabilities of $X(t)$, 
satisfying the initial conditions 
$$
 p_{k,n}(0)=\delta_{k,n}
 =\left\{
 \begin{array}{ll}
 1, & k=n, \\
 0, & k\neq n.
 \end{array}
 \right.
$$
Let $Q$ be the infinitesimal generator of the transition function\eq{1}, 
i.e.\ the matrix whose $(k,n)$-th finite entries are:  
\begin{equation}
 q_{k,n}={{\rm d}\over {\rm d}t}p_{k,n}(t)\Big|_{t=0},
 \label{equation:3}
\end{equation}
satisfying the following relations: 
{\em (a)} $q_{k,n}\geq 0$ for all $k,n\in\Stati$ such that $k\neq n$, 
{\em (b)} $q_{n,n}\leq 0$ for all $n\in\Stati$, and  
{\em (c)} $\sum_{n\in\Stati} q_{k,n}=0$ for all $k\in\Stati$. 
\par
The spatial symmetry of Markov processes allows one to approach effectively the 
first-passage-time problem. Indeed, it has been often exploited by various authors 
to obtain closed-form results for first-passage-time distributions; 
see Giorno {\em et al.} \cite{GiNoRi89} and Di~Crescenzo {\em et al.} \cite{DiGiNoRi97} for 
one-dimensional diffusion processes, Di~Crescenzo {\em et al.} \cite{DiGiNoRi95} 
for two-dimensional diffusion processes, and Di~Crescenzo \cite{Di96} for a class 
of two-dimensional random walks. Moreover, in Di~Crescenzo \cite{Di98}  
a symmetry for truncated birth-death processes was expressed 
as in\eq{2}, with $x_i$ suitably depending on the birth and death rates. 
Such symmetry notion can be extended to the wider class of continuous-time 
Markov chains considered above. Indeed, for a set of positive real numbers 
$\{x_n; \;n\in\Stati\}$ there holds:
\begin{equation}
 p_{N-k,N-n}(t)={x_n\over x_k}\,p_{k,n}(t)\qquad 
 \hbox{for all $k,n\in\Stati$ and $t\geq 0$} 
 \label{equation:9}
\end{equation}
if and only if 
\begin{equation}
 q_{N-k,N-n}={x_n\over x_k}\,q_{k,n}\qquad
 \hbox{for all }k,n\in\Stati. 
 \label{equation:10}
\end{equation}
The proof is similar to that of Theorem 2.1 in Di~Crescenzo \cite{Di98}, 
and thus is omitted. 
\par
Eq.\eq{9} focuses on a symmetry with respect to $N/2$, which identifies with the mid point 
of $\Stati$ when $\Stati=\{0,1,\ldots,N\}$. For each sample-path of $X(t)$ from $k$ to $n$ 
there is a symmetric path from $N-k$ to $N-n$, and the ratio of their probabilities is 
time-independent. Hence, in the following we shall say that $X(t)$ possesses a 
{\em central symmetry\/} if relation\eq{9} is satisfied. 
\begin{remark} 
If $X(t)$ possesses a central symmetry, then 
$$
 {x_n\over x_k}={x_{N-k}\over x_{N-n}}\qquad \hbox{for all }k,n\in\Stati.
$$ 
\end{remark}
An example of a Markov chain with finite state-space and a central symmetry 
is given hereafter. 
\begin{example}{\rm
Let $X(t)$ be a continuous-time Markov chain with state-space 
$\Stati=\{0,1,2,3\}$, with $0$ and $3$ absorbing states, and infinitesimal generator   
$$
 Q=\left[
 \begin{array}{cccc}
 0 & 0 & 0 & 0 \cr
 \alpha\,\varrho_0 +\beta & 
 -\alpha(\varrho_0 +\varrho^2)-\beta \,\varrho_0  & 
 \beta\varrho & 
 \alpha\varrho^2\cr
 \alpha & 
 \beta & 
 -\alpha(\varrho_0 +\varrho^2)-\beta\,\varrho_0  & 
 (\alpha\,\varrho_0 +\beta)\varrho \cr
 0 & 0 & 0 & 0 
 \end{array}\right], 
$$
with $\alpha,\beta,\varrho>0$ and $\varrho_0 =1+\varrho$. Then, $X(t)$ has a central 
symmetry, with $p_{N-k,N-n}(t)=\varrho^{k-n}\,p_{k,n}(t)$ for all $k,n\in\Stati$ and 
$t\geq 0$, and $q_{N-k,N-n}(t)=\varrho^{k-n}\,q_{k,n}$ for all $k,n\in\Stati$. 
}\end{example}
\begin{remark}\label{remark1}
If $X(t)$ has a central symmetry and possesses a stationary distribution 
$\{\pi_n,\;n\in\Stati\}$, with $\ds\lim_{t\to +\infty}p_{k,n}(t)=\pi_n>0$ 
for all $k,n\in\Stati$, then the following statements hold: \\
(a) \ Sequence $\{x_n\}$ is constant, so that $p_{N-k,N-n}(t)=p_{k,n}(t)$ for all 
$k,n\in\Stati$ and $t\geq 0$. \\
(b) \ The stationary distribution is symmetric with respect to $N/2$, i.e. 
$$
 \pi_{N-n}=\pi_n\qquad 
 \hbox{for all $n\in\Stati$.} 
$$
(c) \ Let $X^*(t)$ be the reversed process of $X(t)$, obtained from 
$X(t)$ when time is reversed, and characterized by rates and transition probabilities 
$$
 q_{k,n}^*={\pi_n\over \pi_k}\,q_{n,k},  
 \qquad 
 p_{k,n}^*(t)={\pi_n\over \pi_k}\,p_{n,k}(t),  
 \qquad k,n\in\Stati,\quad t\geq 0.
$$
Then, also $X^*(t)$ has a central symmetry, with 
$p_{N-k,N-n}^*(t)=p_{k,n}^*(t)$ for all $k,n\in\Stati$ and $t\geq 0$. \\ 
(d) \ Let $D=\{d_{k,n}\}$ be the {\em deviation matrix\/} of $X(t)$, with 
elements (see Coolen-Schrijner and Van Doorn \cite{Co02}) 
$$
 d_{k,n}=\int_0^{+\infty}[p_{k,n}(t)-\pi_n]\,{\rm d}t, 
 \qquad k,n\in\Stati.
$$
Then, $D$ has a central symmetry, i.e.\ $d_{N-k,N-n}=d_{k,n}$ for all $k,n\in\Stati$. 
\end{remark}
\par
An example of a Markov chain satisfying the assumptions of Remark~\ref{remark1} is 
the birth-death process on $\Stati$ with birth rate $\lambda_n=\alpha\,(N-n)$ 
and death rate $\mu_n=\alpha\,n$ (see Giorno {\em et al.} \cite{GiNeNo85}, 
or Section 4.1 of Di~Crescenzo \cite{Di98}). 
\section{First-passage-time densities} 
In this section we shall focus on the first-passage-time problem for Markov chains $X(t)$ 
that have a central symmetry and that satisfy the following assumptions:
\newline
(i) \ $N=2s$, with $s$ a positive integer; 
\newline
(ii) \ $q_{i,j}=q_{j,i}=0$, $\sum_{i\in \Stati_-}q_{i,s}>0$, $\sum_{j\in \Stati_+} q_{j,s}>0$, 
$\sum_{i\in \Stati_-}q_{s,i}>0$ and $\sum_{j\in \Stati_+}q_{s,j}>0$ for all $i\in \Stati_-$ 
and $j\in \Stati_+$, where 
$$
 S_-=\{n\in\Stati;\; n<s\}, \qquad 
 S_+=\{n\in\Stati;\; n>s\};
$$
(in other words, if states $i$ and $j$ are separed by  $s$ then 
all sample-paths of $X(t)$ from $i$ to $j$, or from $j$ to $i$, must cross $s$); 
\newline
(iii) \ the subchains defined on $S_-$ and $S_+$ are irreducibles.  
\par
In addition, we introduce the following non-negative 
random variables:
$$
 \begin{array}{l}
 T^{+}_{i,s}= 
 \hbox{ upward first-passage time of $X(t)$ from state $i\in S_-$ to state $s$}, \\
 T^{-}_{j,s}=
 \hbox{ downward first-passage time of $X(t)$ from state $j\in S_+$ to state $s$}. 
 \end{array}
$$
We shall denote by $g^+_{i,s}(t)$ and $g^-_{j,s}(t)$ the corresponding 
probability density functions. Due to assumptions (i)-(iii), for all $t>0$ 
such densities satisfy the following renewal equations:
\begin{eqnarray} 
 && p_{i,j}(t)= \int_0^t g^+_{i,s}(\vartheta)\,p_{s,j}(t-\vartheta)\,{\rm d}\vartheta, 
 \qquad i\in S_-, \;\; j\in \{s\}\cup S_+,
 \label{equation:22} \\
 && p_{j,i}(t)= \int_0^t g^-_{j,s}(\vartheta)\,p_{s,i}(t-\vartheta)\,{\rm d}\vartheta,
 \qquad i\in S_-\cup\{s\}, \;\; j\in S_+.
 \label{equation:23}
\end{eqnarray}
For all $t>0$ and $k\in\Stati$ let us now introduce the {\em probability currents\/}  
\begin{eqnarray}
 && \hspace{-1cm} h^+_{k,s}(t)
 =\lim_{\tau\downarrow 0}{1\over\tau}\,\Prob\{X(t+\tau)=s,\,X(t)<s\mid X(0)=k\} 
 =\sum_{i\in \Stati_-} p_{k,i}(t)\,q_{i,s},
 \label{equation:28} \\
 && \hspace{-1cm} h^-_{k,s}(t)
 =\lim_{\tau\downarrow 0}{1\over\tau}\,\Prob\{X(t+\tau)=s,\,X(t)>s\mid X(0)=k\}
 =\sum_{j\in \Stati_+} p_{k,j}(t)\,q_{j,s}.
 \label{equation:29} 
\end{eqnarray}
They represent respectively the upward and downward entrance probability fluxes at state $s$ 
at time $t$. Due to assumptions (i)-(iii) and Eqs.\eq{22}-(\ref{equation:29}), for 
$i\in \Stati_-$, $j\in \Stati_+$ and $t>0$ they satisfy the following integral equations:
\begin{eqnarray} 
 && \hspace{-1cm} h_{i,s}^-(t)=
\int_0^t g^+_{i,s}(\vartheta)\,h_{s,s}^-(t-\vartheta)\,{\rm d}\vartheta,
 \label{equation:30}\\
 && \hspace{-1cm}
 h_{j,s}^+(t)=\int_0^t g^-_{j,s}(\vartheta)\,h_{s,s}^+(t-\vartheta)\,{\rm d}\vartheta.
 \label{equation:43}
\end{eqnarray}
\par
Hereafter we extend Proposition 2.2 of Di~Crescenzo \cite{Di98} to the case of Markov chains. 
\begin{proposition}\label{proposition1}
Under assumptions (i)-(iii), for all $i\in \Stati_-$, $j\in \Stati_+$ and $t>0$ 
the following equations hold:
\begin{eqnarray}
 && g^+_{i,s}(t)
 =h^+_{i,s}(t)-\int_0^t g^+_{i,s}(\vartheta)\,h_{s,s}^+(t-\vartheta)\,{\rm d}\vartheta,
 \label{equation:24} \\
 && g^-_{j,s} (t)
 =h^-_{j,s}(t)-\int_0^t g^-_{j,s}(\vartheta)\,h_{s,s}^-(t-\vartheta)\,{\rm d}\vartheta.
 \label{equation:25}
\end{eqnarray}
\end{proposition}
\begin{proof}
For all $t>0$ and $i\in \Stati_-$, making use of assumptions (i)-(iii) and Eq.\eq{28} we have
$$ 
 {{\rm d}\over {\rm d}t} p_{i,s}(t)
 =\sum_{n\in\Stati} p_{i,n}(t)\,q_{n,s} 
 =h^+_{i,s}(t)+\sum_{n\in\{s\}\cup S_+} p_{i,n}(t)\,q_{n,s}.
$$
Hence, recalling\eq{22} we obtain
\begin{eqnarray}
 && \hspace{-1.1cm}
 h^+_{i,s}(t)={{\rm d}\over {\rm d}t}
 \left[\int_0^t g^+_{i,s}(\vartheta)\,p_{s,s}(t-\vartheta)\,{\rm d}\vartheta\right]
 -\sum_{n\in\{s\}\cup S_+}\left[\int_0^t g^+_{i,s}(\vartheta)\,
 p_{s,n}(t-\vartheta)\,{\rm d}\vartheta\right]q_{n,s} 
 \nonumber \\
 && =g^+_{i,s}(t)+\int_0^t g^+_{i,s}(\vartheta)
 \left[\frac{\partial}{\partial t}p_{s,s}(t-\vartheta) 
 -\sum_{n\in\{s\}\cup S_+}p_{s,n}(t-\vartheta)\,q_{n,s}\right]{\rm d}\vartheta,
 \label{equation:21}
\end{eqnarray}
where use of initial condition $p_{s,s}(0)=1$ has been made. From Chapman-Kolmogorov 
forward equation we have
$$
 {\partial\over\partial t}p_{s,s}(t-\vartheta) 
 -\sum_{n\in\{s\}\cup S_+}p_{s,n}(t-\vartheta)\,q_{n,s}
 =h^+_{s,s}(t-\vartheta),
 \qquad t>\vartheta, 
$$
so that Eq.\eq{21} gives\eq{24}. The proof of\eq{25} goes along similar lines. 
\end{proof}
\par 
With reference to a Markov chain that has a central symmetry, we now come to the main 
result of this paper, expressing the first-passage-time densities through the symmetry 
state $s$ as difference of probability currents\eq{28} and\eq{29}. 
\begin{theorem}\label{theorem3}
For a Markov chain that has a central symmetry and satisfies assumptions (i)-(iii), 
for all $t>0$ and $k\in \Stati$ there results:
\begin{equation}
 h^-_{2s-k,s}(t)=\frac{x_s}{x_k}\,h^+_{k,s}(t).
 \label{equation:27}
\end{equation}
Moreover, for all $i\in \Stati_-$, $j\in \Stati_+$ and $t>0$ the 
upward and downward first-passage-time densities through state $s$ are given by 
\begin{equation}
 g^+_{i,s}(t)=h_{i,s}^+(t)-h_{i,s}^-(t),
 \qquad 
 g^-_{j,s}(t)=h_{j,s}^-(t)-h_{j,s}^+(t).
 \label{equation:26}
\end{equation}
\end{theorem}
\begin{proof}
Recalling that $N=2s$, for $t>0$ we have 
\begin{eqnarray*}
 && h^-_{2s-k,s}(t)=\sum_{j\in\Stati_+} p_{2s-k,j}(t)\,q_{j,s}
 \hspace{2.1cm} \hbox{(from\eq{29})} \\
 && \hspace{1.7cm}
 =\sum_{i\in\Stati_-} p_{2s-k,2s-i}(t)\,q_{2s-i,s}
 \hspace{1.1cm} \hbox{(setting $j=2s-i$)} \\
 && \hspace{1.7cm}
 ={x_s\over x_k}\sum_{i\in\Stati_-} p_{k,i}(t)\,q_{i,s}
 \hspace{2.1cm} \hbox{(from\eq{9} and\eq{10})} \\
 && \hspace{1.7cm}
 ={x_s\over x_k}\,h^+_{k,s}(t).
 \hspace{3.4cm} \hbox{(from\eq{28})}
\end{eqnarray*}
Eq.\eq{27} then holds. In particular, for $k=s$ it implies that 
$h^-_{s,s}(t-\vartheta)=h^+_{s,s}(t-\vartheta)$ for all $t>\vartheta$. 
Hence, relations\eq{26} follow from Eqs.\eq{30}-(\ref{equation:25}). 
\end{proof}
\par
For a Markov chain $X(t)$ satisfying assumptions (i)-(iii) let us now 
introduce the {\em $s$-avoiding transition probabilities\/}: 
$$
 p^{\langle s\rangle}_{k,n}(t)
 =\Prob\left\{X(t)=n,\,X(\vartheta)\neq s \,\hbox{ for all } \vartheta \in(0,t) 
 \mid X(0)=k\right\},
$$
where $k,n\in\Stati_-\cup\Stati_+$. We note that $p^{\langle s\rangle}_{k,n}(t)$ 
is related to $p_{k,n}(t)$ by 
\begin{equation}
 p^{\langle s\rangle}_{k,n}(t)
 =\left\{
 \begin{array}{ll}
 p_{k,n}(t)-\ds\int_0^t g^+_{k,s}(\vartheta)p_{s,n}(t-\vartheta)\,{\rm d}\vartheta, 
 &  \quad k,n\in\Stati_-, \\
 \hfill & \hfill \\
 p_{k,n}(t)-\ds\int_0^t g^-_{k,s}(\vartheta)p_{s,n}(t-\vartheta)\,{\rm d}\vartheta, 
 &  \quad k,n\in\Stati_+.
 \end{array}
 \right.
 \label{equation:31}
\end{equation}
In the following theorem, for symmetric Markov chains two different expressions 
are given for $p^{\langle s\rangle}_{k,n}(t)$ in terms of $p_{k,n}(t)$. 
It extends Theorem 2.4 of Di~Crescenzo \cite{Di98}; 
the proof is similar and therefore is omitted. 
\begin{theorem}\label{theorem4}
Under the assumptions of Theorem \ref{theorem3}, for  $t>0$ and for 
$k,n\in\Stati_-\cup\Stati_+$ there holds:
\begin{eqnarray*} 
 && p^{\langle s\rangle}_{k,n}(t) = p_{k,n}(t)-\ds{x_k \over x_s}\,p_{2s-k,n}(t)  \\
 && \hspace{1.15cm} 
 = p_{k,n}(t)-\ds{x_s \over x_n}\,p_{k,2s-n}(t).
\end{eqnarray*}
\end{theorem}
\par
We conclude this section by pointing out that for a Markov chain having a central 
symmetry, for all $t>0$ the following relations hold:
\begin{equation}
 g^+_{i,s}(t)=\frac{x_i}{x_s}g^-_{2s-i,s}(t),
 \qquad  
 g^-_{j,s}(t)=\frac{x_j}{x_s}g^+_{2s-j,s}(t),
 \qquad i\in \Stati_-, \;\; j\in \Stati_+, 
 \label{equation:36}
\end{equation}
$$
 p^{\langle s\rangle}_{2s-k,2s-n}(t)=\frac{x_n}{x_k}\,p^{\langle s\rangle}_{k,n}(t),
 \qquad k,n\in\Stati_-\cup\Stati_+.
$$
\section{A bilateral birth-death process with jumps} 
In this section we shall apply the above results to a special symmetric Markov chain 
$X(t)$ with state-space $\Interi$, characterized by the 
following transitions: (a) from $n\in\Interi$ to $n+1$ with rate $\lambda$, 
(b) from $n\in\Interi$ to $n-1$ with rate $\mu$, and (c) from $n\in\Interi-\{0\}$ 
to $0$ with rate $\alpha$. Hence, $X(t)$ is a bilateral birth-death process 
that includes jumps toward state $0$. In order to obtain an expression for the 
transition probabilities $p_{k,n}(t)$, we note that for all $t>0$ the following system holds:
\begin{eqnarray*}  
 && {{\rm d}\over {\rm d}t} p_{k,n}(t)  
 = -(\lambda+\mu+\alpha)\,p_{k,n}(t)+\lambda\,p_{k,n-1}(t)+\mu\,p_{k,n+1}(t), 
 \qquad n\in\Interi-\{0\}, \\
 && {{\rm d}\over {\rm d}t} p_{k,0}(t)
 = -(\lambda+\mu)\,p_{k,0}(t)+\lambda\,p_{k,-1}(t)+\mu\,p_{k,1}(t)
 +\alpha\sum_{r\neq 0} p_{k,r}(t).
\end{eqnarray*}
The probability generating function 
$$
 H(z,t)=\sum_{n=-\infty}^{+\infty} p_{k,n}(t)\,z^n
$$
is thus solution of 
\begin{equation}
 {\partial \over \partial t} H(z,t) =u(z)\,H(z,t)+\alpha,
 \label{equation:35}
\end{equation}
where $u(z)=-(\lambda+\mu+\alpha)+\lambda\,z+\ds\frac{\mu}{z}$, 
with initial condition $H(z,0)=z^k$. 
The unique solution of\eq{35} is 
\begin{equation}
 H(z,t)=H(z,0)\,e^{u(z)\,t}+ \alpha \int_0^t e^{u(z)\,\tau}\,{\rm d}\tau.
 \label{equation:32}
\end{equation}
Hence, recalling that 
$$
 \exp \left\{\left(\lambda\,z+\frac{\mu}{z}\right)t\right\}
 =\sum_{n=-\infty}^{+\infty} I_n(\gamma\,t)\,(\beta\,z)^n  
$$
for $\gamma=2 \sqrt{\lambda\,\mu}$ and $\beta=\sqrt{{\lambda}/{\mu}}$, 
from~(\ref{equation:32}) we obtain:
\begin{equation}
 H(s,t)
 =e^{ -(\lambda+\mu+\alpha)t}\sum_{n=-\infty}^{+\infty} I_{n-k}(\gamma t)\,\beta^{n-k}\,s^n 
 +\alpha \int_0^t e^{ -(\lambda+\mu+\alpha)\tau}  
 \sum_{n=-\infty}^{+\infty} I_n(\gamma \tau)\,(\beta s)^n \,{\rm d}\tau,
\label{equation:44}
\end{equation}
where $I_n(x)$ denotes the modified Bessel function of the first kind. 
Equating the coefficients of $z^n$ on both sides of\eq{44} finally 
yields the transition probabilities
\begin{equation}
 p_{k,n}(t)
 ={\left(\frac{\lambda}{\mu}\right)}^{\!\!\frac{n-k}{2}}
 I_{n-k}\left(2 \sqrt{\lambda \mu} \,t\right) \, e^{-(\lambda+\mu+\alpha)t}
 + \alpha {\left(\frac{\lambda}{\mu}\right)}^{\!\!\frac{n}{2}} \!
 \int_0^t e^{ -(\lambda+\mu+\alpha)\tau}   
 I_n( 2 \sqrt{\lambda \mu} \,\tau) \,{\rm d}\tau.
 \label{equation:33}
\end{equation}
Note that\eq{33} can be expressed as 
\begin{equation}
 p_{k,n}(t)= e^{-\alpha t}\,\widehat{p}_{k,n}(t)
 + \alpha \int_0^t e^{ -\alpha \tau} \,\widehat{p}_{0,n}(\tau)\,{\rm d}\tau,
 \label{equation:34}
\end{equation}
where, for all $t\geq 0$ and $k,n\in\Interi$, 
\begin{equation}
 \widehat{p}_{k,n}(t):={\left(\frac{\lambda}{\mu}\right)}^{\!\!\frac{n-k}{2}}
 I_{n-k}(2 \sqrt{\lambda \mu} \,t) \,e^{-(\lambda+\mu)t}, 
 \label{equation:37}
\end{equation}
is the transition probability of the Poisson bilateral birth-death process with 
birth rate $\lambda$ and death rate $\mu$ (see, for instance, Section~2.1 of Conolly \cite{Co75}). 
Assuming that the stationary probabilities $\pi_n=\ds\lim_{t \to +\infty} p_{k,n}(t)$ 
exist for all $n\in \Interi$, from\eq{35} we have 
\begin{eqnarray*}
 && \sum_{n=-\infty}^{+\infty}\pi_n\,z^n 
 =\lim_{t\to +\infty} H(z,t)
 =-{\alpha \over u(z)}
 ={\alpha z \over \lambda (z-z_1) (z_2-z)} \\
 && \hspace{2cm}
 ={\alpha \over \lambda (z_2-z_1)}
 \Bigg[\sum_{n=-\infty}^{-1}\Big({z\over z_1}\Big)^{\! n}
 +\sum_{n=0}^{+\infty}\Big({z\over z_2}\Big)^{\! n}\Bigg], 
\end{eqnarray*}
where  
$$
 z_{1,2}={\lambda+\mu+\alpha\pm\sqrt{(\lambda+\mu+\alpha)^2-4\lambda\mu}\over 2\lambda},
 \qquad 0<z_1<1<z_2.
$$
Hence,
\begin{equation}
 \pi_n=\left\{
 \begin{array}{ll}
 \ds{\alpha \,{z_1}^{-n}\over {\lambda(z_2-z_1)}}  & \hbox{ for } n=-1,-2,\ldots, \\
 \hfill & \hfill \\
 \ds{\alpha \,{z_2}^{-n}\over {\lambda(z_2-z_1)}}  & \hbox{ for } n=0,1,2,\ldots.
 \end{array} 
 \right. 
 \label{equation:40}
\end{equation}
\par 
It is not hard to see that if $\lambda=\mu$ then $X(t)$ has a central symmetry 
with respect to state $0$, with $x_k=1$ for all $k$: 
$$
 p_{-k,-n}(t)=p_{k,n}(t), \qquad q_{-k,-n}=q_{k,n} 
$$
for all $t>0$ and $k,n\in \Interi$. Note that if $\lambda=\mu$, then $z_1$ and $z_2$ are 
reciprocal zeroes of $u(z)$, so that the stationary distribution\eq{40} is symmetric, 
i.e.\ $\pi_n=\pi_{-n}$ for all $n\in\Interi$. Since $q_{i,j}=q_{j,i}=0$, 
$q_{i,0}>0$ and $q_{j,0}>0$ for all $i,j\in\Interi$ such that $i<0<j$, and 
$q_{0,-1}>0$ and $q_{0,1}>0$, this Markov chain satisfies assumptions (i)-(iii) for which $0$ 
is a symmetry state. In this case the 
first-passage-time densities through $0$ can be obtained via Theorem~\ref{theorem3}. 
Indeed, if $\lambda=\mu$, making use of\eq{33} 
and of property $I_n(x)=I_{-n}(x)$, for all $t>0$ and $k=1,2,\ldots$ we have: 
\begin{eqnarray}
 && \hspace{-1.1cm} 
 g^-_{k,0}(t)= h_{k,0}^-(t)-h_{k,0}^+(t)
 =\sum_{j=1}^{+\infty} p_{k,j}(t)\,q_{j,0}-\sum_{i=-\infty}^{-1} p_{k,i}(t)\,q_{i,0} 
 \nonumber \\
 && = \lambda \left[ p_{k,1}(t) - p_{k,-1}(t) \right]
 + \alpha \Bigg[\sum_{j=1}^{+\infty}p_{k,j}(t) - \sum_{i=-\infty}^{-1}p_{k,i}(t) \Bigg] 
 \nonumber \\
 && =e^{- (2\lambda+\alpha) t} 
 \Bigg\{\lambda\,\left[I_{k-1}(2\lambda\,t)-I_{k+1}(2\lambda\,t)\right]  
 + \alpha \sum_{j=1}^{+\infty} \left[I_{k-j}(2\lambda\,t) - I_{k+j}(2\lambda\,t) \right]\Bigg\}. 
 \label{equation:41} 
\end{eqnarray}
Furthermore, recalling\eq{36}, in this special case for all $t>0$ and $k=1,2,\ldots$ there holds:
$$
 g^+_{-k,0}(t)=g^-_{k,0}(t).
$$  
In analogy with Theorem \ref{theorem4} and by virtue of\eq{33}, when $\lambda=\mu$,  we have
\begin{eqnarray}
 && p^{\langle 0\rangle}_{k,n}(t) 
 = p_{k,n}(t)-p_{-k,n}(t)  
 \nonumber \\
 && \hspace{1.1cm}
 = e^{-(2\lambda+\alpha) t}\,\left[I_{n-k}(2\lambda\,t) -I_{n+k}(2\lambda\,t) \right], 
 \qquad t>0.
 \label{equation:42}
\end{eqnarray}
Note that 
\begin{eqnarray}
&&p^{\langle 0\rangle}_{k,n}(t)=p^{\langle 0\rangle}_{n,k}(t),
\label{equation:45}\\
 &&p^{\langle 0\rangle}_{k,n}(t)
 =e^{-\alpha t}\,\widehat p^{\langle 0\rangle}_{k,n}(t),
 \label{equation:38}
\end{eqnarray}
where $\widehat p^{\langle 0\rangle}_{k,n}(t)$ is the transition probability of 
$\widehat X(t)$ when $\lambda=\mu$. Functions\eq{41} and\eq{42} 
are shown in Figure 1 for some choices of the involved parameters.  
\par
We finally remark that Eqs.\eq{34} and\eq{38} are in agreement with similar 
results for birth-death processes with catastrophes obtained in 
Di~Crescenzo {\em et al.} \cite{DiGiNoRi03a} and \cite{DiGiNoRi03b}. 
\begin{figure}[t]
\begin{center}
\leavevmode
\epsfxsize=350pt
\epsffile{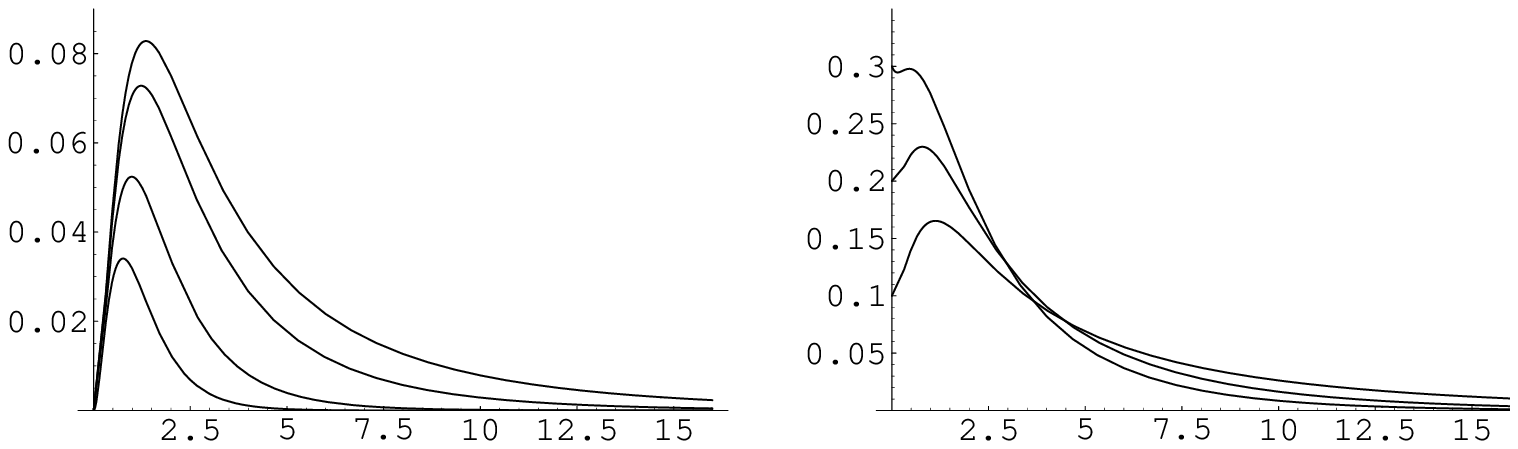}
\caption{On the left-hand are the plots of the downward first-passage-time 
density\eq{41} for $k=3$, $\lambda=1$ and $\alpha=0.1, 0.2, 0.3$, from bottom 
to top near the origin. On the right the $0$-avoiding transition probabilities\eq{42} 
for $k=3$, $n=1$, $\lambda=1$ and $\alpha=0.1, 0.2, 0.5, 1$ (top to bottom) 
are indicated.}
\end{center}
\end{figure}
%
\section{Strong similarity} 
The notion of similarity between stochastic processes has attracted the attention 
of several authors (see Giorno {\em et al.} \cite{GiNoRi88}, for time-homogeneous diffusion 
processes, Guti\'errez J\'aimez {\em et al.} \cite{GuGoRo91} for time-nonhomogeneous 
diffusion processes, Di~Crescenzo \cite{Di94a}, \cite{Di94b}, and Lenin {\em et al.} 
\cite{LePaScVa00}, for birth-death processes, and Pollett \cite{Po01}, for Markov chains). 
Two continuous-time Markov chains $X(t)$ and $\widetilde X(t)$, with state-space 
$\Stati$, are said to be strongly similar if their transition probabilities satisfy
\begin{equation}
 \widetilde p_{k,n}(t)={\beta_n\over \beta_k}\,p_{k,n}(t),
 \qquad \hbox{for all $t\geq 0$ and $k,n\in\Stati$,}
 \label{equation:39}
\end{equation}
where $\{\beta_n,\,n\in\Stati\}$ is a suitable sequence of real positive numbers 
(we refer the reader to Pollett \cite{Po01}, for further details). 
In the following theorem we state that if a Markov chain has a central symmetry, 
then any of its similar chains has a central symmetry as well. 
\begin{theorem}\label{theorem5}
Let $X(t)$ and $\widetilde X(t)$ be strongly similar continuous-time Markov chains 
with state-space $\Stati$; if $X(t)$ has a central symmetry, then for all $t\geq 0$ 
and $k,n\in\Stati$ one has:
$$
 \widetilde p_{N-k,N-n}(t)
 ={\widetilde x_n\over \widetilde x_k}\,\widetilde p_{k,n}(t)
$$
with
$$
 \widetilde x_n={\beta_{N-n}\over \beta_n}\,x_n, \qquad n\in\Stati.
$$
\end{theorem}
The proof is an immediate consequence of assumed symmetry and similarity properties. 
\par
Hereafter we show an application of Theorem \ref{theorem5} to a 
birth-death process having constant rates and state-space $\Interi$. 
\begin{example}\label{esempiob}{\rm
Let $X(t)$ be the bilateral birth-death process with birth and death rates 
$\lambda$ and $\mu$, respectively. From transition probabilities\eq{37} it is 
not hard to see $X(t)$ has a central symmetry with respect to $0$, i.e.\ 
for all $t\geq 0$ and $k,n\in\Interi$ there results
$$
 p_{-k,-n}(t)={x_n\over x_k}\,p_{k,n}(t), 
 \qquad \hbox{with }
 x_n=\left({\lambda\over\mu}\right)^{-n}.
$$
The Markov chains that are strongly similar to $X(t)$ constitute a family 
of bilateral birth-death processes characterized by birth and death rates 
(see Section~4 of Di~Cre\-scenzo \cite{Di94b}, and Example~3 of Pollett \cite{Po01}) 
$$
 \widetilde\lambda_n={\beta_{n+1}\over \beta_n}\,\lambda,
 \qquad
 \widetilde\mu_n={\beta_{n-1}\over \beta_n}\,\mu,
 \qquad n\in\Interi,
$$
and by transition probabilities\eq{39}, with $p_{k,n}(t)$ given in\eq{37} and 
$$
 \beta_n=1+\eta\,\left({\lambda\over\mu}\right)^n, 
 \qquad n\in\Interi,
$$
for all $\eta\geq 0$. Due to Theorem~\ref{theorem5}, the family of strongly 
similar processes has a central symmetry with respect to $0$: 
$$
 \widetilde p_{-k,-n}(t)={\widetilde x_n\over \widetilde x_k}\,\widetilde p_{k,n}(t), 
 \qquad \hbox{with }
 \widetilde x_n={\beta_{-n}\over \beta_n}\,x_n
 ={1+\eta\,\left({\lambda\over\mu}\right)^{-n}\over 1+\eta\,\left({\lambda\over\mu}\right)^n}\, 
 \left({\lambda\over\mu}\right)^{-n},
 \qquad n\in\Interi.
$$
}\end{example}
%

\newpage

\end{document}